\documentclass[12pt]{amsart}
\usepackage{amsaddr}

\usepackage{amsrefs,amssymb,amscd,xcolor}
\usepackage{tikz, subfigure}
\usepackage{mathrsfs}

\usepackage{a4wide}
\usepackage[normalem]{ulem}

\newtheorem{prop}{Proposition}[section]

\newtheorem{theorem}{Theorem}[section]
\newtheorem{lemma}{Lemma}[section]
\newtheorem{cor}{Corollary}[section]
\newtheorem{remark}{Remark}[section]
\theoremstyle{definition}
\newtheorem{defn}{Definition}[section]

\newcommand{\D}{\mathbb{D}}
\newcommand{\R}{\mathbb{R}}
\newcommand{\C}{\mathbb{C}}
\newcommand{\N}{\mathbb{N}}
\renewcommand{\P}{\mathbb{P}}
\newcommand{\E}{\mathbb{E}}
\newcommand{\Z}{\mathbb{Z}}
\newcommand{\I}{{\mathcal I}}
\newcommand{\J}{{\mathcal J}}
\newcommand{\HH}{{\mathcal H}}

\def\L{\mathcal{L}}

\def\LG{\mathcal{L}_G}

\newcommand{\A}{\mathcal A}

\newcommand{\B}{\mathcal B}

\newcommand{\bn}{{\boldsymbol n}}
\def\mybf #1{\textbf{\textit{#1}}}
\title{Invariant densities for random continued fractions
}
\author{Charlene Kalle}
\address{Mathematical Institute, Leiden University, The Netherlands}
\email{kallecccj@math.leidenuniv.nl}
\author{Valentin Matache}
\address{Department of Mathematics, University of Nebraska Omaha, USA}
\email{vmatache@unomaha.edu}
\author{Masato Tsujii}
\address{ Faculty of Mathematics, Kyushu University, Japan}
\email{tsujii@math.kyushu-u.ac.jp}
\author{Evgeny Verbitskiy}
\address{Mathematical Institute, Leiden University, The Netherlands}
\address{ \textit{and} Bernoulli Institute, Groningen University, The Netherlands}
\email{evgeny@math.leidenuniv.nl}

\begin{document}

\begin{abstract} We continue the study of random continued fraction expansions, generated by random application of the Gauss and the R\'enyi backward continued fraction maps. We show that this
random dynamical system admits a unique absolutely continuous invariant measure with smooth density.
\end{abstract}
\date{}
\maketitle

\section{Introduction}
Dynamical systems are traditionally used to generate expansions of real numbers, e.g.,
the so-called $\beta$-expansions
\begin{equation}\label{betaexp}\aligned
x &= \sum_{k\ge 1} \frac  {a_k}{\beta^k},\quad \beta >1, \ a_k\in\{0,1,\ldots, \lfloor \beta\rfloor\},\\
\endaligned
\end{equation}
or the continued fraction expansions
$$\aligned
x &=   \cfrac{1}{a_1
          \pm \cfrac{1}{a_2
          \pm \cfrac{1}{a_3 \pm \cfrac{1}{\ddots} } } }, \quad a_k\in \mathbb N.
\endaligned
$$
The corresponding $\beta$- and continued fractions transformations are classical
objects of study in the theory of dynamical systems. An ergodic point of view (the study of
properties of the invariant measures of these transformations) provides further insights into
the number-theoretic properties of these expansions.

\medskip
In case of $\beta$-expansions, it turned out that for a non-integer base $\beta>1$, Lebesgue almost all numbers in the interval $[0,\frac {\lfloor \beta\rfloor}{\beta-1}]$ admit {\bf uncountably} many different $\beta$-expansions (\cite{EJK,Si}). The natural question became whether one could devise a dynamical way to describe (generate) all possible $\beta$-expansions of a given number $x$. In \cite{Da1} Dajani and Kraaikamp suggested a method based on {\bf random} applications of the
so-called \emph{Greedy} and \emph{Lazy} maps. Indeed, random iterations of these two maps produce all possible $\beta$-expansions of a given number. Further properties of the random  $\beta$-expansions have been established in \cites{Da2,Da3,Da4,Da5,Ke1}.

\medskip

A similar idea can be applied to generate continued fraction expansions of real numbers in $(0,1)$.
In \cite{Ka1} a random continued fraction transformation has been introduced. The two \emph{base} maps are
\begin{itemize}
\item {\it the Gauss continued fraction map}
$$\aligned
T_0(x) &=\left\{\frac 1x\right\},\quad x\in (0,1],
\quad  T_{0}(0)=0,\\
\endaligned
$$
\item {\it the R\'enyi backward  continued fraction map}
$$\aligned
T_1(x) &=\left\{\frac 1{1-x}\right\},\quad x\in [0,1),
\quad T_1(1)=0,
\endaligned
$$
\end{itemize}
where $\{\cdot\}$ denotes the fractional part. We see both maps in Figure~\ref{f:cfmaps} below.

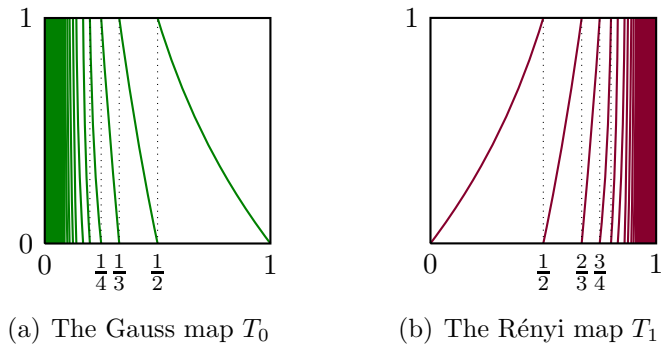
\begin{figure}[h]
\subfigure[The Gauss map $T_0$ ]{
\begin{tikzpicture}[scale=3]
\filldraw[fill=green!50!black, draw=green!50!black] (0,0) rectangle (.09,1);
\draw[thick, green!50!black] (1,0) .. controls (.75,.33) and (.6,.67) .. (.5,1);
\draw[thick, green!50!black] (.5,0) .. controls (.42,.38) and (.36,.78) .. (.33,1);
\draw[thick, green!50!black] (.33,0) .. controls (.29,.45) and (.27,.7) .. (.25,1);
\draw[thick, green!50!black] (.25,0) .. controls (.225,.25) and (.21,.5) .. (.2,1);
\draw[thick, green!50!black] (.2,0) .. controls (.185,.25) and (.175,.5) .. (.17,1);
\draw[thick, green!50!black] (.17,0) .. controls (.155,.25) and (.145,.5) .. (.14,1);
\draw[thick, green!50!black] (.14,0) .. controls (.133,.25) and (.13,.5) .. (.125,1);
\draw[thick, green!50!black] (.125,0) .. controls (.118,.25) and (.115,.5) .. (.11,1);
\draw[thick, green!50!black] (.11,0) .. controls (.105,.25) and (.102,.5) .. (.1,1);
\draw[thick, green!50!black] (.1,0) .. controls (.095,.25) and (.092,.5) .. (.09,1);
\draw[dotted](.5,0)--(.5,1)(.33,0)--(.33,1)(.25,0)--(.25,1)(.2,0)--(.2,1);

\draw[thick](0,0)node[below]{\small 0}--(.25,0)node[below]{\small $\frac14$}--(.33,0)node[below]{\small $\frac13$}--(.5,0)node[below]{\small $\frac12$}--(1,0)node[below]{\small 1}--(1,1)--(0,1)node[left]{\small 1}--(0,0)node[left]{\small 0};
\end{tikzpicture}}
 \hspace{1cm}
\subfigure[The R\'enyi map $T_1$]{
\begin{tikzpicture}[scale=3]
\filldraw[fill=purple!70!black, draw=purple!70!black] (0,-1) rectangle (-.09,0);
\draw[thick, purple!70!black] (-1,-1) .. controls (-.75,-.67) and (-.6,-.33) .. (-.5,0);
\draw[thick, purple!70!black] (-.5,-1) .. controls (-.42,-.62) and (-.36,-.22) .. (-.33,0);
\draw[thick, purple!70!black] (-.33,-1) .. controls (-.29,-.55) and (-.27,-.3) .. (-.25,0);
\draw[thick, purple!70!black] (-.25,-1) .. controls (-.225,-.75) and (-.21,-.5) .. (-.2,0);
\draw[thick, purple!70!black] (-.2,-1) .. controls (-.185,-.75) and (-.175,-.5) .. (-.17,0);
\draw[thick, purple!70!black] (-.17,-1) .. controls (-.155,-.75) and (-.145,-.5) .. (-.14,0);
\draw[thick, purple!70!black] (-.14,-1) .. controls (-.133,-.75) and (-.13,-.5) .. (-.125,0);
\draw[thick, purple!70!black] (-.125,-1) .. controls (-.118,-.75) and (-.115,-.5) .. (-.11,0);
\draw[thick, purple!70!black] (-.11,-1) .. controls (-.105,-.75) and (-.102,-.5) .. (-.1,0);
\draw[thick, purple!70!black] (-.1,-1) .. controls (-.095,-.75) and (-.092,-.5) .. (-.09,0);

\draw[dotted](-.5,-1)--(-.5,0)(-.33,-1)--(-.33,0)(-.25,-1)--(-.25,0)(-.2,-1)--(-.2,0);

\draw[thick](-1,-1)node[below]{\small 0}--(-.5,-1)node[below]{\small $\frac12$}--(-.33,-1)node[below]{\small $\frac23$}--(-.25,-1)node[below]{\small $\frac34$}--(0,-1)node[below]{\small 1}--(0,0)--(-1,0)node[left]{\small 1}--(-1,-1);
\end{tikzpicture}}
\label{f:cfmaps}
\caption{The two continued fraction maps that are used in the random continued fraction system.}
\end{figure}

It is well-known that $T_0$ admits a unique absolutely continuous invariant probability
measure $\mu_0$ with density $\frac{d\mu_0}{d\lambda}(x)=\frac {1}{(1+x)\log 2}$, while $T_1$ only admits a {\bf $\sigma$-finite} absolutely continuous invariant measure $\mu_1$ with density $\frac{d\mu_1}{d\lambda}(x)=\frac 1x$. Here $\lambda$ denotes the one-dimensional Lebesgue measure. The source of \emph{singularity} is the presence of an indifferent fixed point for $T_1$ at the origin.

The random continued fraction map $\mathcal T$ is defined as a skew product transformation from $\{0,1\}^{\mathbb N}\times [0,1]$ into itself, by
\begin{equation}\label{eq:rds}
  \mathcal T(\omega,x) = \big( \sigma\omega, T_{\omega_1}(x) \big),
\end{equation}
where $\sigma:\{0,1\}^{\N}\to\{0,1\}^{\N}$ is the left shift.
For $(\omega,x) \in \{0,1\}^{\mathbb N}\times [0,1]$, write
\[ a_1 = a_1(\omega,x)=m, \, \text{ if } \omega_1 + (-1)^{\omega_1}x\in   \left( \frac{1}{m+1} , \frac1m \right],
\text{ and }  a_k=a_1\big( \mathcal T^{k-1}(\omega,x) \big)\ \text{ for } k>1. \]
Let $\pi: \{0,1\}^{\mathbb N}\times [0,1] \to [0,1]$ be the canonical projection onto the second coordinate. Then for each $n \ge 1$,
\[ x =  \omega_1 + \cfrac{(-1)^{\omega_1}}{a_1 + \omega_2 + \cfrac{(-1)^{\omega_2}}{a_2 + \omega_3 + \ddots + \cfrac{(-1)^{\omega_n}}{a_n+ \pi(\mathcal T^n(\omega,x))}}}. \]

In \cite{Ka1} it is shown that this process converges as $n \to \infty$ and that
\[ x =  \omega_1 +\cfrac{(-1)^{\omega_1}}{a_1 + \omega_2 + \cfrac{(-1)^{\omega_2}}{a_2 + \omega_3 + \ddots }}.
\]

The principal question is whether, for a given $p\in [0,1]$, there exists a $\mathcal T$-invariant
measure $\nu_p$ of the form $m_p\otimes \mu_p$, where $m_p$ is the $(p,1-p)$-Bernoulli measure on $\{0,1\}^{\mathbb N}$, so
\[ m_p([0]) =p,\quad \text{and} \quad m_p([1])=1-p,\]
for the cylinders $[0], [1] \subseteq \{0,1\}^{\mathbb N}$, and $\mu_p$ is an absolutely continuous probability measure on $[0,1]$.
Note that the invariance of $m_p \otimes \mu_p$ is equivalent to the following ``invariance'' condition for $\mu_p$:
\[  \mu_p(A) =p \cdot\mu_p(T_{0}^{-1}A)+(1-p)\cdot \mu_p(T_{1}^{-1}A) \quad\text{ for all Borel sets } A\subseteq [0,1].\]

Clearly, for $p=1$ the answer is positive as the question boils down to the question about the standard Gauss map; similarly, for $p=0$, the answer is negative. The following result has been established in \cite{Ka1}.

\begin{theorem}[Theorem 3.2 and Proposition 3.3 of \cite{Ka1}]\label{res}
For any $p\in (0,1)$ there exists an absolutely continuous invariant measure $\mu_p$ whose density $h_p$ is strictly positive and belongs to the class of functions with {\bf bounded variation}.
\end{theorem}

The density $h_p$ of $\mu_p$ is necessarily a fixed point of the random Perron-Frobenius transfer operator
\begin{equation}\label{maineq}
\L_p f (x) =  \sum_{n=1}^\infty \left[ \frac{p}{(n+x)^2} f \left( \frac1{n+x} \right) + \frac{1-p}{(n+x)^2} f \left( 1-\frac1{n+x} \right) \right],
\end{equation}
which is a weighted average of the Perron-Frobenius operators
\begin{equation}\label{q:gaussrenyipf}
\mathcal L_G f(x) = \sum_{n \ge 1} \frac1{(n+x)^2}  f \left( \frac1{n+x} \right)\quad \text{and} \quad \mathcal L_R f(x) = \sum_{n \ge 1} \frac1{(n+x)^2} f \left( 1- \frac1{(n+x)} \right)
\end{equation}
of the Gauss and R\'enyi maps, respectively. Note that for any $p$, the operator $\L_p$ is {\bf Markov} in the sense that for any $f \in L^1 ([0,1])$,
\begin{itemize}
\item $\L_p f\ge 0$ if $f\ge 0$, and,
\item $\int_{[0,1]} \L_pf(x) \, \lambda(dx) = \int_{[0,1]} f(x) \, \lambda(dx)$.
\end{itemize}
The proof of Theorem~\ref{res} in \cite{Ka1} is based on an application of a theorem by Inoue \cite{Ino}, who established quasi-compactness of transfer operators on spaces of functions of bounded variation for countably branched skew product systems that are expanding on average. In \cite{Ka1} it was conjectured that the invariant probability density $h_p$ is in fact  smooth.

{ 
In the present paper we will study the properties of $\L_p$ on different spaces of smooth functions, namely
 \begin{itemize}
\item
 the space of $k$-times continuously differentiable functions on $[0,1]$,
 \item the Banach space of bounded analytic functions on a certain disk $\mathbb D\subset \C$,
 \item the  Hardy space of analytic functions
on the half-plane $\C_+=\{\textsf{Re}\, z>0\}$.
\end{itemize} Our first result is the following. 

\begin{theorem}\label{goodbound} For any $p\in (0,1)$   the Gauss-R\'enyi transfer operator $\L_p$, given by (\ref{maineq}), is a well defined bounded linear operator on $C^k([0,1])$ for any $k\ge 1$.
Moreover, the essential spectral radius of $\L_p$ on $C^k([0,1])$ satisfies
$$
r_\text{ess}\Bigl( \L_p|_{C^k}\Bigr) \le \zeta(2k+2)-\min(p,1-p),
$$
where $\zeta$ is the Riemann zeta-function.
\end{theorem}

The proof of this result relies on Theorem \ref{ck}, which provides an upper bound for the essential spectral radius of transfer-type operators of the form
\begin{equation}\label{gentrans}
\mathcal L f(x) =\sum_{n\in\I} a_{n}(x) f(b_n(x)),
\end{equation}
where $\I$ is an at most countable index set and $a_n:[0,1] \to \mathbb R$, $b_n :[0,1] \to [0,1]$ for all $n$. Namely, the essential spectral radius of $\L|_{C^k([0,1])}$ satisfies
$$
r_{ess}(\L|_{C^k([0,1])}) \le \limsup_{m\to\infty} \left(
\sup_{x\in[0,1]} \sum_{\bn\in\I^m}
|a_{\bn}(x)|\cdot|b'_{\bn}(x)|^k\right) ^{1/m},
$$
where $\{ a_{\bn}, b_{\bn}:\,\bn\in\I^m\}$ are the `coefficients' of $\L^m$:
$\L^mf=\sum_{\bn\in\I^m}
a_{\bn}\cdot f\circ b_{\bn}$, $m\ge 1$.

The bound provided by Theorem \ref{ck} is not novel, e.g., \cite{BaBook, Isola,Nus16}.
However, contrary to the previous works we do not assume that all maps
$b_n:[0,1]\to [0,1]$ in (\ref{gentrans}) are strict contractions.

 As an immediate corollary of Theorem \ref{goodbound}, we are able to conclude that $\L_p$ is quasi-compact
on spaces $C^k([0,1])$, since $r_{ess}<1$ for $k$  large enough.
Taking into account the results from \cite{Ka1} on uniqueness of  absolutely continuous invariant measures, we, therefore, can conclude that the invariant density  $h_p\in C^k([0,1])$ for all $k$.

}

At first sight it might seem that the random dynamical system built using  the `good' Gauss map and the `bad' R\'enyi map, has an invariant density which is as smooth
as the invariant density of the Gauss map. However, this is not the case.
There is an actual loss of `smoothness' due to the presence of the R\'enyi map: for example, $h_p$ is not real-analytic. We will explain in greater detail
how the presence of the indifferent fixed point of the R\'enyi map
affects the smoothness of the invariant density.

In order to study  properties of $\L_p$ on spaces of analytic functions, we
will employ the technique of  \emph{modification of Markov operators}.
Namely, we will represent the Markov operator $\L_p$, as the sum of two non-negative sub-Markov operators
\[ \mathcal L_p f(x) = \A_p f(x)+\B_pf(x),\]
where $\B_p$ contains some non-hyperbolic (not strictly contracting) inverse branches, and we will consider the modified transfer operator
\[ \widehat{\mathcal L_p} f(x) =\A_p (1-\B_p)^{-1} =\sum_{m=0}^\infty \A_p\B_p^mf(x).\]
Provided $\widehat{\mathcal L_p}$ is well-defined, it is not very difficult to show that $\widehat{\mathcal L_p}$ is again a Markov operator. If we are able to find a positive invariant density $\widehat h_p$ for $\widehat{\mathcal L_p}$, then
\[ h_p=\sum_{m=0}^\infty \B_p^m \widehat h_p\]
can be shown to be an invariant density for $\L_p$ (Proposition \ref{prop:mod}).

In Section \ref{s:banach} we apply this method to study the operator $\L_p$ on the Banach space $\HH^\infty(\mathbb D)$ of analytic bounded functions on the disk $\mathbb D$ which has the interval $[0,1]$ as its diameter.
We isolate the first  branches of the Gauss and R\'enyi maps by setting
\[ \B_p f(z)=  \frac{p}{(1+z)^2} f \left( \frac1{1+z} \right) + \frac{1-p}{(1+z)^2} f \left( 1-\frac1{1+z}\right), \quad \A_p f(z)=\L_p f(z)-\B_p f(z).
\]

\begin{theorem}\label{thm:banachcomplex} For every $p\in (0,1)$ the following two properties hold.
\begin{itemize}
\item[(i)] The operator $\A_p$ is {\bf nuclear} on $\HH^\infty(\mathbb D)$.
\item[(ii)] The operator $J_p=(1-\B_p)^{-1}$ is bounded on $\HH^\infty(\mathbb D)$.
\end{itemize}

\end{theorem}

Again, as an immediate corollary of Theorem \ref{thm:banachcomplex},
one concludes that the operator $\widehat\L_p=\A_p J_p$ is then also nuclear, which easily allows us to conclude that any invariant probability density $\widehat h_p$ of $\widehat\L_p$, and hence, $h_p$ of $\L_p$, are elements of $\HH^\infty(\mathbb D)$.

Finally, in Section \ref{s:hilbert} we turn to the Hardy space $H^2(\C_+)$ of analytic functions on the right half-plane $\C_+=\{z:\textsf{Re} z>0\}$.
We will consider
\[ \B_p f(z)=  \frac{1-p}{(1+z)^2} f \left( 1-\frac1{1+z}\right), \quad \A_p f(z)=\L_p f(z)-\B_p f(z).
\]

\begin{theorem} For every $p\in (0,1)$ the following two properties hold.
\begin{itemize}
\item[(i)] $\A_p$ is a compact operator on  $H^2(\C_+)$.
\item[(ii)] $\widehat{\L_p}=\A_p(1-\B_p)^{-1}$ are {\bf compact} operators
on the Hardy space $H^2(\C_+)$.
\end{itemize}
\end{theorem}

Similarly to the previous cases, as an immediate corollary we can conclude that that
the eigenfunction $\widehat h_p$ of the compact Markov operator $\widehat{\L_p}=\A_p(1-\B_p)^{-1}$, and hence also $h_p$ of $\L_p$, can be extended to an analytic function in the Hardy space $H^2(\C_+)$.

Throughout the article we use $\mathcal L_p$ to denote the operator from \eqref{maineq} related to the random continued fraction system and $\mathcal L$ will denote an operator in general.

\subsection*{Acknowledgments.} Research of CK was partially supported by NWO Veni-grant number 639.031.140. CK and EV are grateful to Kyushu University and the World Premier International Researcher Invitation Program “Progress 100” for hospitality and support. EV is grateful to A.J.E.M. Janssen, D. Terhesiu for helpful discussions.

 \section{Spectral Properties of Transfer Operators,
Sub-Markovian Operators, and Jump Transformations}
\label{s:overview}

Suppose $H$ is a Banach space, and  $K\subset X$ is a convex cone such that $K-K$ is dense in $H$. An operator $\L:H\to H$ is called {\bf positive} (with respect to the cone $K$)
if $\L K\subseteq K$.
The celebrated Krein–Rutman theorem states
that if $\L$ is a compact positive operator
with a positive spectral radius $r(\L)>0$,
then $r(\L)$ is an eigenvalue of $\L$, and there exists a `positive' $h\in K$, $h\ne 0$, such that $\L h=r(\L) h$.

Transfer operators of  dynamical systems are rarely compact.
There is a substantial amount of literature devoted to the problem of existence of invariant densities of transfer operators, see e.g.,  the book by Baladi \cite{BaBook}
and references therein. One of the most popular approaches is based on establishing
quasi-compactness of transfer operators.

\begin{defn}An operator $\L$ acting on a Banach space $H$ is called \mybf{quasi-compact} if there exists an $r<r(\L)$
such that the Banach space $H$ can be decomposed as $
H = G\oplus F$,
where $G$ and $F$ are $\L$-invariant closed  subspaces such that
$$
\dim G<\infty \text{ and }  r( \L |_F)<r,
$$
or, equivalently, $\L$ has a finite number (counting with multiplicity)
of eigenvalues of absolute value $\ge  r$.
The essential spectral radius $r_{ess}(\L)$ of $\L$ is defined as the greatest lower bound of such $r$'s.
\end{defn}

Nussbaum extended the Krein-Rutman theorem to quasi-compact positive operators in \cite{Nussbaum2}.
For the Banach spaces considered in this paper, namely $C^k([0,1])$ and the Hardy space $H^2(\C_+)$,
the positive cone $K$ is simply
the set of all non-negative functions on $[0,1]$:
$$
K=\{f\in H: f(x)\ge 0\quad \forall x\in [0,1]\}.
$$

Therefore, if we are able to show that our operator $\L_p$ is compact or quasi-compact,
we immediately conclude that the spectral radius $r(\L_p)$ is an eigenvalue with a non-negative eigenfunction. Since $\L_p$ is also Markov, we necessarily have that $r(\L_p)=1$,
since $\int_{[0,1]} \L_p h(x) \, \lambda(dx)=\int_{[0,1]} h(x) \, \lambda(dx)$ for all non-negative $h$.

A principal question is of course whether $r(\L)$ is a simple eigenvalue of the positive operator $\L$. For a general positive operator,  \emph{strong positivity} or \emph{strict quasi-positivity} of $\L$ is sufficient \cite{Sasser}. For transfer operators $\L$ of dynamical systems some form of mixing of the underlying dynamics is sufficient to ensure
simplicity of the maximal eigenvalue. Indeed, our random continued fraction dynamical system possesses the necessary mixing properties, and in \cite{Ka1} we showed that it  has  a unique absolutely continuous invariant measure $m_p\otimes \mu_p$, and hence, necessarily, $\L_p$ has a unique invariant probability density $h_p$.

In the present paper we will establish compactness of $\L_p$ on two Banach spaces of analytic functions.
On $C^k([0,1])$, $k\ge 1$ sufficiently large, $\L_p$ will be shown to be quasi-compact.

There are a number of standard methods to establish quasi-compactness of transfer-like operators. The most popular way  is to establish a Doeblin-Fortet or Lasota-Yorke type inequality \cite{HennionBook}:
if for some $n\ge 1$, $r\in (0,r(\L))$, $C>0$, and all $f\in H$ one has
$$
\|\L^nf \|\le r^n \|f\|+C\|f\|_w,
$$
where $\|\cdot \|_w$ is a weaker (semi-)norm on $H$ such that the unit ball in the strong norm is relatively compact in the weak norm, then $\L$ is quasi-compact. The quasi-compactness of $\L_p$ has been established by such means on spaces $C^1([0,1])$ and $C^2([0,1])$ for $p$  sufficiently close to $1$ (i.e., $1-p$ sufficiently small) in \cite{BRS}, {  and recently extended in \cite{TC} to all $C^k([0,1])$, $k\ge 2$.}

Another  approach is based on the well-known formula of Nussbaum \cite{Nussbaum} for the essential spectral radius: if $\L$ is a bounded linear operator
on a Banach space $H$, then
$$
r_{ess}(\L) = \lim_{t\to\infty}\left(
\inf\Big\{ \|\L^t -\mathcal K\| \, \big|  \, \mathcal K:H\to H\text{ is compact}\Big\} \right)^{1/t}.
$$
Thus, if one is able to show that $r_{ess}(\L)<r(\L)$, then $\L$ is quasi-compact.

Building on the formula for the spectral radius of transfer operators of smooth expanding interval maps obtained by  Collet \& Isola \cite{CI}, we will derive a similar  upper bound on the essential spectral radius of rather general transfer-type operators of the form
$\L f(x) =\sum_{n\in\mathcal I} a_n(x) f(b_n(x))$ acting on $C^k$. However, we will
not make any assumptions on the contraction rates of the $b_n$'s.

Applying  our bound to $\L_p$, $p\in (0,1)$, acting on $C^k([0,1])$, we will show
that
\[ r_{ess}(\L_p|_{C^k([0,1])})\le \zeta(2k+2)-\min(p,1-p),\]
where $\zeta$ is the Riemann zeta-function. Since for $p\in (0,1)$, $\min(p,1-p)>0$,  one concludes that $r_{ess}(\L_p|_{C^k([0,1])})<1$ for all sufficiently large $k$, and hence, $\L_p$ has a  positive invariant probability density $h_p\in C^k([0,1])$
for all sufficiently large, and thus for all, $k\ge 1$.

\subsection{Analytic approach}

If one is interested in further spectral and  analytic properties of transfer operators and their invariant  densities, it is often useful to consider the action of transfer-type operators
\begin{equation}\label{analyt}
\L f(z) =\sum_{n=1}^\infty a_n(z) f( b_n(z)),\ z\in\Omega,
\end{equation}
on analytic functions on a certain domain $\Omega\subset \C$. Here, $a_n:\Omega\to\C$, $b_n:\Omega\to\Omega$ are assumed to be analytic. Operators as in \eqref{analyt} are sums of \emph{weighted composition operators}. Ruelle \cite{Rue76} has observed that if the $b_n$'s are contractions, then the corresponding transfer operator is \emph{nuclear} on the appropriate space of analytic functions, see \cite{MR2} for the treatment of the Gauss map and \cite{BJ} for recent rather general results and a comprehensive overview.

Babenko has proposed a novel approach to the study of the transfer operator of the Gauss map
\cite{Bab} acting on certain Hilbert spaces of analytic functions on half-planes.
The approach was further developed for Gauss and Gauss-type transfer operators \cite{MR1,Mayer, Jen1} and related maps \cite{Isola}.

We apply both methods to the study of the appropriate Markovian modifications of $\L_p$, which we introduce now, in Sections \ref{s:banach} and \ref{s:hilbert}.

\subsection{Submarkov operators}

The main technical difficulty in the analysis of our transfer operator $\L_p$
stems from the presence of two indifferent (non-expanding) points: $x=1$ for the Gauss map,
and $x=0$ for the R\'enyi map, see Figure \ref{f:cfmaps}.
The point $x=1$ is not a fixed point of $T_0$,
thus one typically considers the second power of $\L_0 = \mathcal L_G$.
The point $x=0$, on the other hand, is an indifferent fixed point of $T_1$.
The standard approach in such situations is to consider \emph{induced} systems,
which often have better hyperbolic properties, and to draw conclusions
about the original system from the corresponding properties of induced
systems, e.g., continue an absolutely continuous invariant probability measure
of the induced system to an absolutely continuous invariant measure
of the original system.

In the present paper we will use a different, although somewhat related method, based on a modification of the Markov operator we want to understand.

\begin{prop}\label{prop:mod} Suppose $(X,\mathcal F,\mu)$ is a probability space, and
$H$ is some Banach space of real-valued functions on $X$ such that $H\subset L^1(X,\mu)$.
Suppose also that $\A$ and $\B$ are non-negative bounded linear operators on $H$ such that
\begin{itemize}
\item their sum $\L=\A+\B$ is a Markov operator, i.e.,
$$
\int_X \L u(x)\mu(dx) = \int_X u(x)\mu(dx) \quad \forall u\ge 0, u\in H,
$$
\item $J=(1-\B)^{-1}=\sum_{m=0}^\infty \B^m$ is a bounded operator on $H$.
\end{itemize}
Then
\begin{equation}\label{transjump}
\widehat \L = \A J =\sum_{m=0}^\infty \A\B^m
\end{equation}
is Markov. Moreover, if $\hat h\in H$, $\hat h\ge 0$, is such that $\widehat \L \, \hat h=\hat h$, then $ h= J \hat h$
is non-negative and satisfies
$$\L h = h.
$$
\end{prop}

\begin{proof}Suppose  $u\ge 0$, then
$\int_{X} \L u (x) \mu(dx)= \int_{X} u(x) \mu(dx)$. On the other hand,
$$\aligned
\int_X \A u(x) \mu(dx)&= \int_X(\L-\B) u(x) \mu(dx) = \int_X \L u(x) \mu(dx )- \int_X \B u(x)\mu(dx)\\
&= \int_{X} u(x) \mu(dx) -\int_X \B u(x) \mu(dx)=\int_{X} (1-\B) u(x) \mu(dx).
\endaligned
$$
Suppose now  $v\ge 0$. Let $u=J v =(1-\B)^{-1}v$, then $u\ge 0$, and applying the previous equality, we get
$$\aligned
\int_{X} {\widehat \L} \, v d\mu=\int_{X} \A (J v) d\mu=
\int_{X} (1-\B) (1-\B)^{-1} v d\mu=\int_X v d\mu,
\endaligned
$$
which means that $\widehat{\L}$ is Markov.
Now suppose that
$\widehat \L \, \hat h=\hat h$ and put $h=J\hat h=\sum_{m=0}^\infty\B^m \hat h$. Then
\[
\L h=\A h+\B h =\A J\hat h +\sum_{m=1}^\infty \B^m \hat h=
\widehat{\L}\, \hat h +\sum_{m=1}^\infty \B^m \hat h=
\hat h +\sum_{m=1}^\infty \B^m \hat h= h. \qedhere
\]
\end{proof}

Thus if we represent a Markov transfer-type operator $\L f(x) =\sum_{n\in I} a_n(x)f (b_n(x))$,
as a sum of two sub-Markov operators
$$
\A f(x) =\sum_{n\in I_0} a_n(x)f (b_n(x)),\quad
\B f(x) =\sum_{n\in I\setminus I_0} a_n(x)f (b_n(x)),
$$
such that for some Banach space $H$
\begin{itemize}
\item $\A$ has good spectral properties, say, $\A$ is compact (or nuclear) on  $H$, and
\item $(1-\B)^{-1}$ is a bounded operator on $H$,
\end{itemize}
then the Markov operator $\widehat \L=\A(1-\B)^{-1}$ has equally good
spectral properties, as a composition of  compact (nuclear) and  bounded operators. And as Proposition \ref{prop:mod} shows, any invariant density $\hat h$ of
$\widehat\L$ gives rise to an invariant density  $h$ of $\L$ in the same Banach space $H$.

In  this paper we will consider two splits of $\L_p$ of
such nature:
$$
\B_p f(x)= \frac{p}{(x+1)^2} f \left( \frac1{x+1} \right) + \frac{1-p}{(x+1)^2} f \left( 1-\frac1{x+1} \right),\quad \A_p f(x)= \L_pf(x)-\B_pf(x)
$$
in Section \ref{s:banach}, and
$$
\B_p f(x)=  \frac{1-p}{(x+1)^2} f \left( 1-\frac1{x+1} \right),\quad \A_p f(x)= \L_pf(x)-\B_pf(x)
$$
in Section \ref{s:hilbert}. The first split `isolates' the two non-uniformly expanding
branches ($n=1$) of the Gauss and R\'enyi maps, while the second split removes
only the problematic branch of the R\'enyi map.

\subsection{Jump transformation}
A natural question is how the proposed method of studying $\widehat\L=\A(1-\B)^{-1}$ compares with the more standard techniques of considering the
induced transformations. As we will demonstrate
now, the two are closely related. The advantage of the proposed method is that one does not necessarily have to understand
the combinatorial aspect of the induced transformation.

Let us illustrate the basic idea with an example \cite[Example 3.4]{Bon}.
Suppose $T:[0,1]\to [0,1]$ is a $C^1$-piecewise expanding map of the interval with two intervals of monotonicity
$\{I_0,I_1\}$ and full branches, i.e.,  such that $TI_0=TI_1=[0,1]$. Consider the first hitting time of $I_1$, $\tau_{I_1}(x)=\inf \{n\ge 0:\ T^nx\in I_1\}$ and let $\widehat T(x)=T^{\tau_{I_1}(x)+1}$ be the jump transformation. Then $\widehat T$ is again a piecewise monotonic map on $[0,1]$, and the transfer operator $\widehat \L$ of $\widehat T$ satisfies, c.f. (\ref{transjump}),
$$
\widehat \L f(x) = \sum_{n=0}^\infty \A\, \B^n f(x),
$$
where
$$
\A f(x)=\frac 1{|T'(x_1)|} f(x_1),\quad \B f(x)=\frac 1{|T'(x_0)|} f(x_0),
$$
 are `transfer operators' corresponding to the two branches of  $T$: here for $x\in [0,1]$,  $x_0=T^{-1}x\cap I_0$ and
$x_1=T^{-1}x\cap I_1$ are the two preimages, and clearly
$$
\L f(x) = \sum_{y\in T^{-1}x} \frac 1{|T'(y)|} f(y)=\A f(x)+\B f(x).
$$
More generally,  for a  measure preserving dynamical system
$(X, \mathcal B, \mu, T)$ and a measurable set $E$ satisfying $T(E)=X$ and $\bigcup_{k \ge 0}T^{-k}(E)=X$,
the jump transformation $T_E: X \to X$ is defined by $T_E(x) = T^{p(x)}(x)$, where the first passage time $p:X \to \mathbb  N \cup \{ \infty \}$ to $E$ is
\[ p(x) = 1+\inf \{ n \ge 0 \, : \, T^n(x) \in E \}.\]

Isola \cite{Isola} studied the Farey map
$$
T(x) =\begin{cases} \frac {x}{1-x},&\ \text{if }\ 0\le x\le \frac 12,\\
\frac {1-x}{x},&\ \text{if }\ \frac 12< x\le 1,\\
\end{cases}
$$
with the transfer operator $\L$ for the Farey map $$
\L f(x)=\left(\frac{1}{1+x}\right)^{2}\left[f\left(\frac{x}{1+x}\right)+f\left(\frac{1}{1+x}\right)\right]=\left(\A+\B\right) f(x).
$$
If $E=[\frac 12,1]$, the jump transformation $T_E$ is the Gauss map.
In order to  investigate the spectra of transfer operators of the Farey and Gauss maps, Isola considered operators $\widehat{\L}_{z}=z \A\left(1-z \B\right)^{-1}$, $z\in\C$, and studied
their properties on  certain Hilbert spaces of holomorphic functions.
For somewhat similar ideas see also \cite{Rugh}.

Proposition \ref{prop:mod} is a rather general result, and in principle, does not require the consideration of any induced or jump transformations. However, the application of this result to the analysis of $\L_p$, $p\in(0,1)$,
does have a somewhat `hidden' inducing mechanism.

\section{Spectral gap and quasi-compactness on $C^{k}([0,1])$}\label{sec:spectralgap}

Consider the following rather general transfer-like operator
\begin{equation}\label{linop}
\mathcal L f (x) =\sum_{n\in \I}
a_n(x) f( b_n(x)),
\end{equation}
where $\I$ is some finite or countable set of indices, and $a_n:[0,1]\to \R$ and $b_n:[0,1]\to [0,1]$ for all $n\in\I$.

For $m\ge 1$ any $\bn =(n_1,\ldots, n_m)\in\I^m$, put
\begin{equation}\label{iterate}
\aligned
 b_{\bn}(x) & = (b_{n_1} \circ \cdots \circ b_{n_m})(x),\\
 a_{\bn}(x) &=
 a_{n_m}(x)\cdot
 a_{n_{m-1}}(b_{n_m}(x))\cdot
 a_{n_{m-2}}(b_{(n_{m-1},n_m)}(x))\cdot \ldots \cdot a_{n_1}( b_{(n_2,\ldots, n_m)}(x))\\
 &= \prod_{i=1}^m (a_{n_i} \circ b_{n_{i+1}} \circ \cdots \circ b_{n_m})(x).
 \endaligned
\end{equation}
Then
$$
\mathcal L^m f(x) = \sum_{\bn\in\I^m} a_{\bn}(x) f(b_{\bn}(x)).
$$

\medskip
Let $C^k([0,1])$ denote the Banach space of all $k$-times
continuously differentiable functions on $[0,1]$ with the norm
$$
\| f \|_k = \max_{j=0,\ldots,k} \sup_{x\in[0,1]} |f^{(j)}(x)|= \max_{j=0,\ldots,k} \| f^{(j)} \|_0.
$$
We now turn to estimating the  essential spectral radius of $\L$ on $C^k([0,1])$.
 
\begin{theorem}\label{ck} Assume that $a_n,b_n\in C^k([0,1])$ for all $n\in \I$,  are such that
\begin{equation}\label{normcondition}
\sum_{n\in\I} \| a_n \|_k \Bigl(1+ \| b_n \|_k\Bigr)^{k}<\infty,
\end{equation}
then $\L$, given by \eqref{linop},
is a bounded linear operator on $C^{k}([0,1])$. Moreover, the essential spectral
radius of $\L|_{C^k([0,1])}$ satisfies
$$
r_{ess}(\L|_{C^k([0,1])}) \le \limsup_{m\to\infty} \left(
\sup_{x\in[0,1]} \sum_{\bn\in\I^m}
|a_{\bn}(x)|\cdot|b'_{\bn}(x)|^k\right) ^{1/m}.
$$
\end{theorem}

{ 
\begin{remark}
The estimate is given by the same formula as in case of
transfer operators of expanding maps \cite{ BaBook, CI, Nus16}.
Note however, that we do not require all maps $b_n:[0,1]\to [0,1]$ to be contractions, which is a standard assumption in previously
published results.
Hence, the theorem in the above form becomes applicable to non-uniformly expanding (random) interval maps.
Naturally, the upper bound on the essential spectral radius becomes useful only
if most $b_n$'s are indeed mostly contracting, or as we shall in the following section,
are contracting on average. The proof of Theorem \ref{ck} is a direct adaptation of the method of Collet and Isola \cite{CI} and is provided for completeness.
\end{remark}
}
\begin{proof} Recall that if $a,f$ are two $k$-times continuously differentiable functions on
$[0,1]$, then
for every $m\in\{0,1,\ldots,k\}$ one has
$$
(a\cdot f)^{(m)}(x) =\sum_{j=0}^m {m\choose j} a^{(j)}(x) f^{(m-j)}(x).
$$
Furthermore, if $b:[0,1]\to [0,1]$ is a $k$-times differentiable function, then for every $m=1,\ldots,k$, by the Fa\`a di Bruno  formula \cite{Faa} one has
$$\aligned
(f\circ b)^{(m)}(x) &= \sum_{j_1+2j_2+\cdots +mj_m=m}
\frac {m!}{j_1! j_2! \cdots j_m!} f^{(j_1+j_2+\cdots+j_m)}( b(x)) \prod_{i=1}^m \left( \frac{b^{(i)}(x)}{i!} \right)^{j_i}\\ 
&=
\sum_{j=1}^m f^{(j)}( b(x)) 
\sum_{\substack{j_1+2j_2+\cdots +mj_m=m,\\ j_1+j_2+\cdots+j_m=j}}
\frac {m!}{j_1! j_2! \cdots j_m!} \prod_{i=1}^m \left( \frac{b^{(i)}(x)}{i!} \right)^{j_i}\\
&=\sum_{j=1}^m f^{(j)}( b(x)) B_{m,j}\left(b'(x),b''(x),\dots,b^{(m-j+1)}(x)\right),
\endaligned
$$
where $B_{m,j}(z_1,\ldots,z_{m-j+1})$ are the Bell polynomials
\[ \begin{split}
B_{m,j}(z_1& ,\ldots, z_{m-j+1})\\
=\ & \sum_{\substack{\ell_1+2\ell_2+\ldots +(m-j+1)\ell_{m-j+1}=m\\ \ell_1+\ell_2+\ldots+\ell_{m-j+1}=j}}
\frac {m!}{\ell_1! \ell_2! \cdots \ell_m!}
\Bigl( \frac {z_1}{1!}\Bigr)^{\ell_1}
\Bigl( \frac {z_2}{2!}\Bigr)^{\ell_2}\cdots
\Bigl( \frac {z_{m-j+1}}{(m-j+1)!}\Bigr)^{\ell_{m-j+1}}.
\end{split}\]

The equalities above imply that
\begin{itemize}
\item  for the product $a\cdot f$ and $m\in\{0,1,\ldots,k\}$ one has
\[ \| (a\cdot f)^{(m)}\|_0 \le 2^{m} \| a \|_{m} \| f \|_{m}\le 2^{k} \| a \|_{k} \| f \|_{k}.\]

\item  for the composition $f\circ b$ and $m=0$ one obviously has
\[  \| f\circ b\|_0\le \| f \|_0,\]
and for $m\in\{1,2,\ldots,k\}$,  one has
\begin{equation}\label{boundFaa}\aligned
\| (f\circ b)^{(m)}\|_0&\le \|f \|_{m} \sum_{j=1}^m B_{m,j}(\| b \|_m,\ldots,\|b \|_m)
\le   \| f \|_k   \sum_{j=1}^m B_{m,j}(1,1,\ldots,1)\| b \|_m^j.
\endaligned
\end{equation}
\end{itemize}

The  Stirling numbers of the second kind are defined as
$$
S_{m,j}=B_{m,j}(1,1,\ldots,1) 
=\frac{1}{j!}\sum_{i=0}^{j}(-1)^{j-i}{j \choose i} i^m,\ j=1,\ldots,m, \quad S_{m,0}=0,
$$
and have the following upper bound:
$$
S_{m,j}\le \frac 12{m\choose j} j^{m-j}, \ j=1,\ldots,m.
$$
Hence, for $x\ge 0$, and $m=1,\ldots,k$,  one has
$$
\sum_{j=1}^m S_{m,j}x^j \le \sum_{j=1}^m \frac 12{m\choose j} j^{m-j} x^j\le
\frac {m^m}2\sum_{j=1}^m  {m\choose j}  x^j\le
\frac {m^m}{2}(1+x)^m\le  {k^k}(1+x)^k.
$$
Thus for every $n\in\I$, by combining the estimates above,  if $m=1,\ldots,k$, one gets
$$\aligned
\left\| \bigl(a_n\cdot f \circ b_n\bigr)^{(m)}\right\|_0&\le \sum_{j=0}^m {m \choose j} \| a_n^{(m-j)}\|_0 \| (f \circ b_n)^{(j)}\|_0\\
&= \| a_{n}^{(m)}\|_0 \cdot\|f\|_0+\sum_{j=1}^m {m\choose j} \| a_{n}^{(m-j)}\|_0 \cdot
\|(f\circ b_n)^{(j)}\|_0\\
&\le \| a_{n} \|_k \cdot\|f\|_k +\sum_{j=1}^m {m\choose j} \| a_{n}\|_k k^k \|f\|_k (1+\|b_n\|_k)^k\\
&\le   k^k \| a_{n} \|_k (1+\|b_n\|_k)^k \left( 1+\sum_{j=1}^m {m\choose j}\right)
\|f\|_k \\
&\le(2k)^k \| a_{n} \|_k  (1+\|b_n\|_k)^k \|f\|_k.
\endaligned
$$
Therefore, if
$$
\sum_{n\in\I} \| a_n \|_k \Bigl(1+ \| b_n \|_k\Bigr)^{k}<\infty,
$$
then
$ \mathcal L f (x) =\sum_{n\in \I}
a_n(x) f( b_n(x))
$
is indeed a bounded linear operator on $C^{k}([0,1])$.
For $t\ge 1$,
$$
\L^t f(x) = \sum_{\bn\in\I^t} a_{\bn}(x) f(b_{\bn}(x)),
$$
where $a_{\bn}, b_{\bn}$ are given by \eqref{iterate}, and thus $\L^t$ is again an operator of the form (\ref{linop}) with
$$
\sum_{\bn\in\I^t} \| a_{\bn}\|_k (1+\|b_\bn\|_k)^k<\infty.
$$
 
We now turn to the estimation of the essential spectral radius. Suppose $\L$ is operator of the form (\ref{linop})
satisfying  the norm condition \eqref{normcondition}
(i.e., $\L$ could be some power of the original operator). We would like to estimate from above
$$
\inf\{\|\L-\mathcal K\|_{k}\, \big| \, \mathcal K:C^k([0,1])\to C^k([0,1])\text{ is compact}\}.
$$

Suppose the  operator $\L$ is of the form $\mathcal L f (x) =\sum_{\bn\in \J} a_\bn(x) f( b_\bn(x))$, then, using expressions for the derivatives of products and compositions of functions, for $m\ge 0$, one has
$$\aligned
(\L  f)^{(m)}(x) &=\Bigl( \sum_{\bn\in\J } a_\bn(x) f(b_\bn(x))\Bigr)^{(m)}=\sum_{\bn\in\J} \sum_{\ell=0}^m {m\choose \ell}a_\bn^{(m-\ell)}(x) (f(b_\bn(x))^{(\ell)}\\
&=\sum_{\bn\in\J} \sum_{j=0}^m f^{(j)}(b_\bn(x))
\Bigl[\sum_{\ell=j}^m {m\choose \ell} a_\bn^{(m-\ell)}(x) B_{\ell,j}\bigl(b_\bn^{(1)}(x),
\ldots, b_\bn^{(\ell-j+1)}(x)\bigr)\Bigr],
\endaligned
$$
where we set $B_{0,0}\equiv 1$, and $B_{\ell,0}\equiv 0$ for $\ell\ge 1$.
Note also, that since $B_{k,k}(z)= z^k$, one has
$$\aligned
(\L  f)^{(k)}(x)&=\sum_{\bn\in\J} a_{\bn}{(x)} \bigl(b_\bn'(x)\bigr)^k f^{(k)}(b_\bn(x))+
\sum_{j=0}^{k-1}\mathcal G_{j,m}  f^{(j)}(x),
\endaligned
$$
where the operators $\mathcal G_{j,m}^{(t)}$, $j=0,\ldots, m-1$, are given by
$$
\mathcal G_{j,m} f^{(j)} (x)=\sum_{\bn\in\J}
\Bigl[\sum_{\ell=j}^m {m\choose \ell} a_{\bn}^{(m-\ell)}(x) B_{\ell,j}\bigl(b_\bn^{(1)}(x), \ldots, b_\bn^{(\ell-j+1)}(x))\Bigr] f^{(j)}(b_\bn(x)).
$$

From this point our proof is a straightforward  adaptation of the proof of \cite[Theorem 2.5]{BaBook}.
We will use the same bijection 
$$
  (C^k([0,1]),\|\cdot \|_k) \ni f\mapsto (f,f',\ldots, f^{(k)})\in (\tilde B,\|\cdot\|_{\tilde B}),
$$
where
\[ \tilde B=\{ \psi=(\psi_0,\ldots,\psi_k):\, \psi_k\in C([0,1]),\ \psi'_m=\psi_{m+1},\ m=0,\ldots,k-1\},\]
and $ \|\psi\|_{\tilde B}=\max_{m=0,\ldots,k} \|\psi_m\|_{0}$.

We will introduce the following notation: if  $\phi_0=\L\psi_0$,
then  $\phi_m=\phi_0^{(m)} = (\L\psi_0)^{(m)}$ and
$$
\begin{bmatrix}
\phi_0\\
\phi_1\\
\vdots\\
\phi_{k-1}\\
\phi_k
\end{bmatrix}=
\begin{bmatrix}
\mathcal G_{0,0} & 0 & 0 &\ldots & 0 & 0\\
\mathcal G_{0,1} & \mathcal G_{1,1} & 0  &\ldots & 0 &0\\
\vdots &\vdots &\vdots & \ddots &\vdots &\vdots\\
\mathcal G_{0,k-1} & \mathcal G_{1,k-1} & \mathcal G_{2,k-1}   &\ldots & \mathcal G_{k-1,k-1}
& 0\\
\mathcal G_{0,k} & \mathcal G_{1,k} & \mathcal G_{2,k}   &\ldots &  \mathcal G_{k-1,k} &\mathcal G_{k,k}\\
\end{bmatrix}
\begin{bmatrix}
\psi_0\\
\psi_1\\
\vdots\\
\psi_{k-1}\\
\psi_k
\end{bmatrix}
=\begin{bmatrix}
Q_0\psi_0\\
Q_1\psi_0\\
\vdots\\
Q_{k-1}\psi_0\\
Q_k\psi_0 +\mathcal G_{k,k} \psi_0^{(k)}
\end{bmatrix}.
$$

The  key observation is that the operators $Q_m$, $m=0,\ldots, k$,
viewed as operators from $C^k([0,1])$ to $C([0,1])$, are compact.  To prove this claim we have to show that the image of the unit sphere in $C^k([0,1])$ is relatively compact in $C([0,1])$.
By the Arzel\'a-Ascoli theorem it is sufficient to check the equicontinuity; boundness is clear.

Hence, for any operator $\L$ of the form
$$\mathcal L f (x) =\sum_{\bn\in \J} a_\bn(x) f( b_\bn(x)),
$$
where $a_\bn,b_\bn$ satisfy the norm condition (\ref{normcondition}), one has
$$\aligned
\inf_{\mathcal K:C^k([0,1])\to C^k([0,1])} \|\L f-\mathcal Kf \|_{k}&\le \left\| \mathcal G_{k,k} f^{(k)}\right\|_0=\sup_{x} \sum_{\bn\in \J}
|a_\bn(x)||b_\bn'(x)|^k |f^{(k)}(x)|\\
&\le \left( \sup_{x} \sum_{\bn\in \J}
a_\bn(x)|b_\bn'(x)|^k\right) \|f\|_k.
\endaligned
$$
Now, applying this bound to the powers $\L^{t} f(x)=\sum_{\bn\in\I^t} a_{\bn}(x)f(b_{\bn}(x))$, i.e., $\J=\I^t$,
we obtain the desired result.
\end{proof}

\subsection{Essential spectral radius of the random Gauss-R\'enyi transfer operator}
Now we are ready to apply Theorem \ref{ck} to the operator
\[ \L_p f(x)=\sum_{n=1}^\infty \Big[ \frac{p}{(n+x)^2} f \Big( \frac1{n+x} \Big) + \frac{1-p}{(n+x)^2} f \Big( 1-\frac1{n+x} \Big) \Big].
\]
This operator can be represented as follows: let
\[ \Omega=\{0,1\},\quad \mathcal I =\mathbb N\times \Omega = \{ (n,\omega):
 \ n\in\N,\omega=0,1\},\]
 and put
\begin{equation}\label{setup}
\aligned
a_{n,\omega}(x)&=\begin{cases}
\frac {p}{(x+n)^2},&\quad\text{if } \omega=0,\\
\frac {1-p}{(x+n)^2},&\quad\text{if } \omega=1,
\end{cases},\quad
b_{n,\omega}(x)=\begin{cases}
\frac {1}{x+n},&\quad\text{if } \omega=0,\\
1-\frac {1}{x+n},&\quad\text{if } \omega=1.
\end{cases}
\endaligned
\end{equation}
Then
\[ \L_p f(x) =\sum_{\tilde n\in\I} a_{\tilde n}(x) f( b_{\tilde n}(x)),\quad \tilde n=(n,\omega).\]
Clearly the functions $a_{\tilde n}$ and $b_{\tilde n}$, $\tilde n=(n,\omega)\in \I$, satisfy the conditions
of Theorem \ref{ck} for all $k\ge 1$. We now turn to estimating the essential spectral radius of $\L_p$.

Note that
\[a_{n,\omega} = \left\{
\begin{array}{ll}
p \cdot |b_{n,\omega}'(x)|, & \text{if } \omega_1=0,\\
(1-p)\cdot |b_{n,\omega}'(x)|, & \text{if } \omega_1=1.
\end{array}\right.\]
Therefore, for every $m\ge 0$ and $k\ge 1$,
$$\aligned
\sum_{\tilde\bn\in\I^m} a_{\tilde \bn}(x) |b_{\tilde \bn}'(x)|^{k} &=\sum_{\bn\in\N^m} \sum_{\omega_1^m\in\{0,1\}^m} \P[\omega_1^m] |b_{\bn,\omega_1^m}'(x)|^{k+1}=\sum_{\bn\in\N^m} \E_\P|b_{\bn,\omega_1^m}'(x)|^{k+1},
\endaligned
$$
where $\P$ is the $(p,1-p)$-Bernoulli measure on $\Omega^{\N}$. Here and below we use the notation $\omega_1^m = (\omega_1, \omega_2, \ldots, \omega_m) \in \{0,1\}^m$ and $n_1^m = (n_1, n_2, \ldots, n_m) \in \mathbb N^m$.

Let us now evaluate the derivative of $b_{\tilde \bn}(x)$, $\tilde\bn=(n_1^m,\omega_1^m)\in\I^m$. First note that for each $(n, \omega ) \in \I$,
\[ b_{n, \omega} (x) = \omega + \frac{(-1)^{\omega}}{x+n} = \frac{\omega x + \omega n + (-1)^{\omega}}{x+n},
\]
so $b_{n, \omega} (x)$ is the M\"obius transformation with matrix $M_{n,\omega}$, where
\[ M_{n,0} = \begin{bmatrix}
0 & 1\\
1 & n
\end{bmatrix} \quad \text{ and } \quad M_{n,1} = \begin{bmatrix}
1 & n-1\\
1 & n
\end{bmatrix},\ n\ge 1.\]
Therefore, $b_{\tilde \bn}=b_{n_1,\omega_1}\circ b_{n_2,\omega_2}\circ \cdots \circ b_{n_m,\omega_m}$ is again a M\"obius transformation, i.e.,
$b_{\tilde \bn}(x)= \frac{Ax +B}{Cx +D}$, where $A,B,C,D \in \mathbb\Z_+$, are the entries of the matrix
\[
\begin{bmatrix} A& B\\
C&D
\end{bmatrix}=M_{n_1, \omega_1} \cdot \ldots \cdot M_{n_{m-1}, \omega_{m-1}}  \cdot M_{n_m, \omega_m}.\]
Note that since $|\det M_{n,0}|=|\det M_{n,1}|=1$ for all $n$, the determinant of the
product is also $\pm 1$, and hence, since $C,D>0$, for the  the derivative one has
\[ \sup_{x \in [0,1]} \big| b_{\tilde \bn}'(x) \big| = \sup_{x\in [0,1]}\Big| \frac{AD-BC}{(Cx+D)^2} \Big| = \frac1{D^2},\]
thus the maximum of all derivatives is attained at $x=0$.

Furthermore, one has
$$
\aligned
|b_{n_1^m,\omega_1^m}'(0)|^{1/2}&=
\cfrac 1{n_m}\cdot
\cfrac 1
{
n_{m-1}+\omega_{m}+
\cfrac {(-1)^{\omega_m}}
{n_m}
}
\times\cdots\\
&\times\cfrac{1}{n_1 +\omega_2 + \cfrac{(-1)^{\omega_2}}{n_2 +\omega_3 +\cfrac{(-1)^{\omega_3} }
{\ddots + \cfrac{(-1)^{\omega_m}}{n_m}}}}
\endaligned
$$
 Hence,
$$
|b_{n_1^m,\omega_1^m}'(0)|=
\Bigl[z_{n_m}(\omega_m)
z_{n_{m-1}^m}(\omega_{m-1}^m)
\cdots z_{n_1^m}(\omega_1^m)\Bigr]^{2},
$$
where
$$
z_k=z_{n_k^m}(\omega_k^m):=
\cfrac{1}{n_k +\omega_{k+1} + \cfrac{(-1)^{\omega_{k+1}}}{n_{k+1} +\omega_{k+2} +\cfrac{(-1)^{\omega_{k+2}} }
{\ddots + \cfrac{(-1)^{\omega_m}}{n_m}}}}.
$$
 
Equivalently, one has
\[ Q_{m}^{(k)} := \sum_{\bn\in\N^m} \E_\P|b_{\bn,\omega_1^m}'(0)|^{k+1}
=\sum_{n_1^m\in\N^m} \mathbb E_{\mathbb P}\Bigl(\Bigl[  z_1 z_2 \cdots z_{m-1}z_m\Bigr]^{2k+2}\Bigr).\]
Let us now estimate
\[ E_j=\mathbb E( z_j^{2k+2}\cdots z_{m}^{2k+2}), \quad j=1,\ldots,m.\]
Let $t=2k+2$. Then by the law of total expectation, one has
$$
\aligned
E_j& = \mathbb E\Bigl( z_j^{t}z_{j+1}^t\cdots z^t_m\Bigr)= \mathbb E\Bigl( \mathbb E\bigl(z_j^{t}z_{j+1}^t\cdots z^t_m \, \bigl|\bigr.\, z_{j+1},\ldots,z_m \bigr)\Bigr)\\
&=\mathbb E\Bigl( \mathbb E\bigl(z_j^{t}\,
\bigl|\bigr.\, z_{j+1},\ldots,z_m) z_{j+1}^t\cdots z^t_m \bigr)\Bigr)=\mathbb E\Bigl( \mathbb E\bigl(z_j^{t}\,
\bigl|\bigr.\,  z_{j+1} )\cdot
z_{j+1}^t\cdots z^t_m \Bigr).
 \endaligned
$$
Furthermore,
$$
\mathbb E\bigl(z_j^{t}\, \bigl|\bigr.\, z_{j+1} )=p\left(\frac{1}{n_j+z_{j+1}}\right)^t+ (1-p)\left(\frac{1}{n_j+1-z_{j+1}}\right)^t
$$
and as a function of $z_{j+1}\in [0,1]$, this expression attains its maximal value at the end points of the interval $[0,1]$. Hence,
\[ \mathbb E\bigl(z_j^t \, \bigl|\bigr.\, z_{j+1} ) \le  \frac {\max(p,1-p)}{n_j^t} +  \frac {\min(p,1-p)}{(n_j+1)^t} =:V_p^{(t)}(n_j).\]
Thus for all $j\in \{1,\ldots, m-1\}$, we get $ E_j\le  V_p^{(t)}(n_j)E_{j+1}$. Hence,
$$
E_1\le V_p^{(t)}(n_1)\cdots V_p^{(t)}(n_{m-1}) \frac 1{n_m^t},
$$
and therefore
$$\aligned
Q_{m}^{(k)}
&=\sum_{n_1^m\in\N^m} \mathbb E_{\mathbb P}\Bigl(\Bigl[
z_{1 }
z_{2}
\cdots z_{m-1}z_m\Bigr]^{2k+2}\Bigr)\\
&\le \sum_{n_1^m\in\N^m} V_p^{(2k+2)}(n_1)\cdots V_p^{(2k+2)}(n_{m-1}) \frac 1{n_m^{2k+2}}=
\left( \sum_{n=1}^\infty V_p^{(2k+2)}(n)\right)^{m-1} \zeta(2k+2),
\endaligned
$$
where $\zeta$ denotes the Riemann zeta-function. Thus  we can conclude that
\[ r_{ess}(\L|_{C^k([0,1])}) \le \limsup_{m\to\infty} \left( \sup_{x\in[0,1]} \sum_{\tilde \bn\in\I^m} a_{\tilde \bn}(x) |b'_{\tilde \bn}(x)|^k\right) ^{1/m}\le \sum_{n=1}^\infty V_p^{(2k+2)}(n).\]
Moreover,
$$\aligned
\sum_{n=1}^\infty V_p^{(2k+2)}(n)&={\max(p,1-p)}\zeta(2k+2)+\min(p,1-p)\bigl(\zeta(2k+2)-1\big)\\
&=\zeta(2k+2) -\min(p,1-p).
\endaligned
$$
For $p\in (0,1)$, $\min(p,1-p)>0$, and since $\zeta(2k+2)\to 1$
as $k\to \infty$, we have that for any $p\in (0,1)$ and all sufficiently
large $k$,
$$
r_{ess}(\L_p|_{C^k([0,1])}) <1.
$$

Therefore, $\L_p$ is a quasi-compact Markov linear operator on $C^k([0,1])$,
and following the arguments in Section \ref{s:overview}, we conclude that $\L_p$ has a
positive fixed point $h_p\in C^k([0,1])$.

\begin{cor} For any $p\in (0,1)$ and $k\ge 1$, the unique invariant probability density $h_p$ of $\L_p$
is $k$-times continuously differentiable.
\end{cor}

\section{The Banach space approach}\label{s:banach}

We have seen that the main difficulty in studying the transfer operator $\L_p$ stems from the fact that the
R\'enyi map has an indifferent fixed point at $x=0$. In this section, following the method described in Section 2, we will consider a certain
modification of the transfer operator $\L_p$. Namely, consider the operators
$$\aligned
\B_pf (x)& = \frac {p}{(1+x)^2} f\left(\frac {1}{1+x}\right)+ \frac {1-p}{(1+x)^2} f\left(1-\frac {1}{1+x}\right), \\
\endaligned
$$
and
$$\aligned
\A_pf (x)& = \sum_{n=2}^\infty \frac {p}{(n+x)^2}
f\left(\frac {1}{n+x}\right)+\sum_{n=2}^\infty  \frac {1-p}{(n+x)^2}  f\left(1-\frac {1}{n+x}\right).
\endaligned
$$

The operator  $\B_p$ isolates the branches of the maps $T_0$ and $T_1$ that are not everywhere expanding.

We will study the behaviour of these operators on a linear Banach space $\mathcal H^\infty(\mathbb D)$ of bounded holomorphic functions on open domains $\mathbb D$, equipped with the
norm
$$
  \| f \| = \sup_{z\in\mathbb D} |f(z)|.
$$

We will denote by $\D[\alpha, \beta]$  the disk in the complex plane $\C$ which has the interval $[\alpha, \beta]\subset\R$, $\alpha<\beta$,  as its diameter.
If the interval $[\alpha, \beta]$ is mapped into $[\alpha', \beta']$ by a M\"obius  transformation $$
\mathbf T(z)=\frac {az+b}{cz+d},\quad a,b,c,d\in\R,
$$ then $\D[a,b]$ is mapped by $\mathbf T$ into $\D[\alpha',\beta']$: $\mathbf T(\D[\alpha,\beta])\subset\D[\alpha',\beta']$.
Moreover, under this condition, if $f\in \mathcal H^\infty(\D[\alpha',\beta'])$, then
$f\circ \mathbf T\in \mathcal H^\infty(\D[\alpha,\beta])$.

If $[\alpha',\beta']\subset[\alpha,\beta]$, then $\D[\alpha',\beta']\subset\D[\alpha,\beta]$,
and hence we have the inclusion
$$
\iota: \mathcal H^\infty(\D[\alpha,\beta]) \to \mathcal H^\infty(\D[a',b']),\quad \iota(f) = f|_{\D[\alpha',\beta']},
$$
and if $\alpha<\alpha'<\beta'<\beta$, i.e., $\D[\alpha',\beta']$ is compactly embedded into
$\D[\alpha,\beta]$, then $\iota$ is \mybf{ nuclear}.
In particular, the inclusion $\iota: \mathcal H (\D [-\frac12,\frac32 ]) \to \mathcal H(\D[0,1])$ is nuclear.

\begin{prop} The operator $\widehat\L_p = \A_p(1- \B_p)^{-1}$ is nuclear on the space $\mathcal H^\infty(\mathbb D[0,1])$.
\end{prop}
\begin{proof} The operators $\ A_p$  and $\B_p$ can be represented as follows
$$\aligned
\A_p f(z)&=p\sum_{n=2}^\infty w_n(z)
f\,\circ \mathbf T_n(z)+
(1-p)\sum_{n=2}^\infty \tilde w_n(z)
f\,\circ  \widetilde{\mathbf T}_n(z)=p\A_p^G f(z)+(1-p)\A_p^Rf(z),\\
 \B_p f(z)&=
p w_1(z)
f\circ \,\mathbf T_1(z)+
(1-p)\tilde w_1(z) \, f\circ \widetilde{\mathbf T}_1(z),
\endaligned
$$
where
$$
\mathbf T_{n}(z) =\frac {1}{z+n},\quad w_n(z)= \frac 1{(z+n)^2},\quad
\widetilde{\mathbf T}_{n}(z) =\frac {-z+n-1}{z+n},\quad
\tilde w_n(z)=\frac 1{(z+n)^2},
\quad n\ge 1.
$$

Clearly, ${\mathbf T}_1(\D[0,1])$ and $\widetilde{\mathbf T}_1(\D[0,1])$ are subsets of
$\D[0,1]$,
and hence $\B_p$ is a bounded operator on $\mathcal H^\infty(\mathbb D[0,1])$.
Moreover, $\B_p^2$ is a contraction. Indeed,
\[ \begin{split}
\B_p^2 f(z) =\, &\frac {p}{(1+z)^2}\Bigl[
\frac {p}{(1+\frac 1{1+z})^2} f\Bigl(\frac {1}{1+\frac 1{1+z}}\Bigr)
+
\frac {1-p}{(1+\frac 1{1+z})^2} f\Bigl(\frac {\frac 1{1+z}}{1+\frac 1{1+z}}\Bigr)
\Bigr]\\
&+\frac {1-p}{(1+z)^2}\Bigl[
\frac {p}{(1+\frac z{1+z})^2} f\Bigl(\frac {1}{1+\frac z{1+z}}\Bigr)
+
\frac {1-p}{(1+\frac z{1+z})^2} f\Bigl(\frac {\frac z{1+z}}{1+\frac z{1+z}}\Bigr)
\Bigr]\\
=\, &\frac {p^2}{(2+z)^2}f\Bigl(\frac {z+1}{z+2}\Bigr)
+\frac {p(1-p)}{(2+z)^2}f\Bigl(\frac {1}{z+2}\Bigr)
+\frac {p(1-p)}{(1+2z)^2}f\Bigl(\frac {1+z}{1+2z}\Bigr)\\
&+\frac { (1-p)^2}{(1+2z)^2}f\Bigl(\frac {z}{1+2z}\Bigr),
\end{split}\]
and hence
$$\|\B_p^2\|\le \frac {p^2}4+\frac {p(1-p)}4+p(1-p)+(1-p)^2=1-\frac 34 p<1.$$
Since $\B_p^2$ is a contraction, the following series converges in the operator norm
$$
  \sum_{n=0}^\infty \B_p^n = (\mathbf 1+\B_p) \Bigl( \sum_{k=0}^\infty \B_p^{2k}\Bigr) =\Bigl( \sum_{k=0}^\infty \B_p^{2k}\Bigr) (\mathbf 1+\B_p)=:J_p
$$
and thus $J_p$ is the inverse of $(1-\B_p)$.

Secondly, every M\"obius transformation ${\mathbf T_n}$ or
$\widetilde{\mathbf T}_n$, $n\ge 2$, maps the interval
$I=[-\frac 12,\frac 32]$ into $[0,1]$:
$$
{\mathbf T_n}I=\left[ \frac {1}{n+\frac 32},\frac {1}{n-\frac 12}\right]\text{ and }
\widetilde{\mathbf T}_nI=\left[ \frac {n-\frac 32}{n-\frac 12},\frac {n+\frac 12}{n+\frac 32}\right].
$$
Thus, if $f\in\mathcal H^{\infty}(\D[0,1])$, then $f\circ{\mathbf T_n}, f\circ\tilde{\mathbf T}_n\in  \mathcal H^{\infty}(\D[-\frac 12,\frac 32])$ for all $n\ge 2$.
Since $w_n,\tilde{w}_n\in \mathcal H^{\infty}(\D[-\frac 12,\frac 32])$ as well, and the series
$$
\sum_{n=2}^\infty w_n(z),\ \sum_{n=2}^\infty \tilde w_n(z)
$$
converge absolutely and uniformly on $\D[-\frac12,\frac32]$, we conclude that
$$\A_p:\mathcal H^{\infty}(\D[0,1])\mapsto \mathcal H^{\infty}(\D[-\frac12,\frac32]).
$$
Therefore, the operator $\A_p$, viewed as an operator from $\HH^\infty(\D[0, 1])$ to $\HH^\infty(\D[0, 1])$, is nuclear,  as
a composition of bounded and nuclear operators:
$$
\HH^\infty(\D[0,1]) \xrightarrow{A_p} \HH^\infty\left(\D\left[-\frac 12,\frac 32\right]\right) \xrightarrow{\iota}\HH^\infty(\D[0,1]).
$$
Finally, $\widehat\L_p= \A_p(1-\B_p)^{-1}:\HH^\infty(\D[0,1])\to\HH^\infty(\D[0,1])$ is also
nuclear -- again, as a composition of nuclear and bounded operators.
\end{proof}

 \begin{cor} The unique invariant probability density $h_p$ on $[0,1]$ of the transfer operator $\L_p$ can be extended to an analytic function in $\HH^\infty(\D[0,1])$.
  \end{cor}

\section{Hilbert space approach}\label{s:hilbert}

The Hilbert space approach to the study of transfer operators of the Gauss and Gauss-type
maps, introduced in \cite{Bab} and developed further in \cite{MR1, Mayer, Isola, Jen1}, consists in identifying an equivalent integral operator acting on an appropriate Hilbert space.
The advantage of the method is that the corresponding integral operator has a continuous symmetric kernel
and hence, a real spectrum. Moreover, since the operator is of trace-class, one can derive
relatively accurate estimates of the second eigenvalue, and hence on the spectral gap.
The method also allows to conclude that the invariant density is analytic on the appropriate
right half-plane in $\C$.

We now summarise the results of the Hilbert space approach to the analysis of the transfer operator of the Gauss map.
Consider the Hilbert space $L^2(\R_+,\mu)$, where $\mu$ is a measure on $\R_+$,
absolutely continuous with respect to the Lebesgue measure, with the density
$$
d\mu(t) =\frac {t}{e^t-1} dt=\frac {te^{-t}}{1-e^{-t}} dt.
$$
The scalar product on $L^2(\R_+,\mu)$ is given by
$$
\langle \phi,\psi\rangle =\int_{0}^\infty   \phi(t)\, \overline{\psi(t)} \, \mu(dt).
$$
Define also the Laplace-Mellin type transform $\widehat \cdot$ on $L^2(\R_+,\mu)$ as follows: for $\phi\in L^2(\R_+,\mu)$ and $x\in [0,1]$, let
$$
{\widehat \phi}(x) =\int_{0}^\infty e^{-tx} \phi(t) \mu(dt).
$$

The following theorem summarises the known results on the Hilbert space method applied to the
 transfer operator of the Gauss map.

\begin{theorem}[Theorem 1 of \cite{MR1}; Lemma 4, 5 and Theorem 1 of \cite{Jen1}] (a) Suppose $\phi\in L^2(\R_+,\mu)$, then $\widehat\phi$ is holomorphic in the right half-plane
$$
R=\Bigl\{z\in\C: \text{\rm Re}(z)>-\frac 12\Bigr\}.
$$
(b) Let $J_1$ denote the Bessel function of the first kind:
$$
J_1(z) =\sum_{k=0}^\infty (-1)^k \frac {\left(\frac z2\right)^{2k+1}}{k!(k+1)!}.
$$
Then the integral operator $\mathcal K$
$$
  \mathcal K\phi(s) =\int_{0}^\infty  \frac {J_1(2\sqrt{st})}{\sqrt{st}}\, \phi(t)\, d\mu(t),\quad s\ge 0,
$$
preserves $L^2(\R_+,\mu)$, is selfadjoint (hence, has real spectrum), and is of trace class.

(c) Moreover, if $\psi=\mathcal K\phi$, with $\phi,\psi\in L^2(\R_+,\mu)$, then $f(x)=\widehat{\phi}(x)$ and
$g(x)=\widehat{\psi}(x)$ satisfy
$$
  g(x) =\sum_{n=1}^\infty \frac 1{(x+n)^2} f\left(\frac 1{x+n}\right),
$$
i.e., $g=\LG f$. Equivalently, the following diagram commutes
$$\begin{CD}
L^2(\R_+,\mu)    @>\mathcal K>>  L^2(\R_+,\mu)\\
@V\widehat{\cdot}VV        @VV\widehat{\cdot}V\\
\mathcal H     @>\LG>>  \mathcal H,
\end{CD}
$$
where $\mathcal H=\widehat{L^2(\R_+,\mu) }$.
\end{theorem}

 \subsection{Hardy spaces and composition operators}

 An operator of type
 $$
   C_\phi f(z) = f\circ\phi(z)=f(\phi(z))\qquad f\in \mathcal{H},
 $$
where $\Omega$ is a nonempty set and $\mathcal{H}$ a space consisting of functions defined on $\Omega$ is called a composition operator. The map $\phi:\Omega\to\Omega$ is called the symbol of $C_\phi $.
A composition operator followed by a multiplication operator is called a
 weighted composition operator. More formally, the operator
 $$
   T_{\psi,\phi} f(z) =\psi(z)\cdot f\circ\phi(z)=\psi(z) f\big(\phi(z)\big)\qquad f\in \mathcal{H},
 $$
is called the weighted composition operator induced by $\psi $ (the weight function, or symbol) and $\phi $ (the composition symbol).

We will be concerned with such operators acting on the Hilbert Hardy spaces over the open unit disc, respectively the right open half-plane.

Let $\mathbb U=\{z \in \mathbb{C}:|z|<1\}$ be   the open unit disk.
We denote by $H^{2}(\mathbb U)$ the space of holomorphic functions $f$ on $\mathbb U$ such that
\[
\|f\|_{H^{2}(\mathbb U)}:=\sup _{0 \leq r<1} \left(\frac{1}{2\pi }\int_{0}^{2\pi}\left|f\left(r e^{i t}\right)\right|^{2} d t\right)^{1 / 2}<\infty.
\]
Equivalently, $H^2(\mathbb U)$ is the space of all functions analytic in the open unit disc, having square--summable Maclaurin coefficients.

The Hardy space $H^2(\C_+)$, $\C_+=\{z:\textsf{ Re }z>0\}$, consists of all analytic functions in the right half-plane, which are square integrable on vertical lines, with respect to Lebesgue measure, with bounded set of integrals:
$$
\|f\|_{H^2(\C_+)} =\sup_{x>0} \left(\frac{1}{\pi }\int _{-\infty }^\infty |f(x+iy)|^2 dy\right)^{\frac 12}<\infty.
$$

\subsection{Norms of weighted composition operators}
A basic, function theory principle, known as  Littlewood's subordination principle \cite{DurenBook}, says that all composition operators on $H^2(\mathbb U)$, induced by symbols fixing the origin, are contractions, that is operators with norm less than or equal to $1$.  Based on this principle, one can establish with little technical effort the norm estimate
 \begin{equation}\label{NormEstimate01}
 \| C_\phi \| \leq \sqrt{\frac{1+|\phi (0)|}{1-|\phi (0)|}},
 \end{equation}
proving that all composition operators induced by analytic selfmaps of the disc are bounded operators acting on $H^2(\mathbb U)$. The case of a half-plane is different. There, only few analytic selfmaps induce bounded composition operators. If one considers weighted composition operators now, the situation becomes more complicated, even when working in $H^2(\mathbb U)$. If any analytic $\psi $ is considered a weight symbol, then note that $T_{\psi,\phi} 1=\psi $ and so, if $\psi $ does not belong to $H^2(\mathbb U)$, then $T_{\psi,\phi} $ is not a bounded operator acting on that space. Furthermore, not all pairs of maps $\psi \in H^2(\mathbb U)$ and $\varphi $ an analytic selfmap of the disc induce bounded, weighted composition operators. General boundedness criteria for such operators do exist, but in complicated terms (such as pull--back Carleson measures), and so, they are hard to use in practical situations. Therefore we resort to the following principle which is sufficient for our needs in this paper. It is well known that, if $\psi $ is a bounded analytic map of the disc, then the multiplication operator induced by it, that is the operator $M_\psi f=\psi f$, is a bounded operator on $H^2(\mathbb U)$ with norm $\| M_\psi \| =\| \psi \| _\infty $. Given the obvious equality $T_{\psi,\phi} =M_\psi C_\phi $, the norm estimate $\| M_\psi C_\phi \| \leq \| M_\psi \| \, \| C_\phi \| $  combines with (\ref{NormEstimate01}) into proving
\begin{equation}\label{NormEstimate02}
\| T_{\psi ,\phi } \| \leq \| \psi \| _\infty \sqrt{\frac{1+|\phi (0)|}{1-|\phi (0)|}},
\end{equation}
{which shows that all weighted composition operators on $H^2(\mathbb U)$, with bounded weight symbol, are bounded.

We will also need an estimate of the  trace norm $\| T\| _1={\rm tr}(\sqrt{T^*T})$, of a weighted composition operators on $H^2(\mathbb U)$.

 \begin{lemma}
If $\psi \in H^2(\mathbb U)$ and $\| \phi \|_\infty <1$, then both $C_\phi $   and
$T_{\psi , \phi }$  are nuclear (or trace class), since the following  estimate holds:
\begin{equation}\label{EquationTrace}
\| T_{\psi , \phi }\| _1\leq \frac{\| \psi \| _{H^2(\mathbb U)}}{1-\| \phi \| _\infty } .
\end{equation}
\end{lemma}
\begin{proof}
Note that
\[
\| T\| _1\leq \sum_{n=0}^\infty \| Te_n\|
\]
where $\{ e_n\} $ is a complete orthonormal basis. Therefore, one has
\[
\| T_{\psi,\phi}\| _1\leq \sum_{n=0}^\infty \| T_{\psi , \phi }(z^n)\|_{H^2(\mathbb U)}=\sum_{n=0}^\infty \| \psi \phi ^n\|_{H^2(\mathbb U)}\leq
\sum_{n=0}^\infty \| \psi \| _{H^{2}(\mathbb U)}\| \phi \| _\infty ^n=\frac{\| \psi \| _{H^{2}(\mathbb U)}}{1-\| \phi \| _\infty } . \qedhere
\]
\end{proof}

{The fact that, if $\psi \in H^2(\mathbb U)$ and $\| \phi \|_\infty <1$, then both $C_\phi $   and
$T_{\psi , \phi }$  are nuclear, appeared originally in \cite[Theorem 2.7]{GKP} with a more complex proof based on results in \cite{HarperSmith} and without the above trace estimate.}

Recently the following result has been established in \cite{Valentin}.

 \begin{prop}\label{propValentin} Let $\phi$ be an analytic selfmap of $\C_{+}$ and $\psi$ an analytic map on the same set. Let $\Phi=\gamma^{-1} \circ \phi \circ \gamma$ be the conformal conjugate of $\phi$ by the Cayley transform $\gamma(z)=\frac {1+z}{1-z} $, $|z|<1$, $\gamma:\mathbb U\to\C_+$,
 and its
inverse $\gamma^{-1}(w) = \frac {w-1}{w+1}$, $\textsf{Re\,} w>0$, $\gamma^{-1}:\C_+\to\mathbb U$.
Denote by $\Psi$ the map
\[
\Psi(z)=\psi \circ \gamma(z) \frac{1-\Phi(z)}{1-z}, \quad z \in \mathbb{U}.
\]
Then the operators $T_{\psi, \phi}:H^2(\C_+)\to H^2(\C_+)$ and $T_{\Psi, \Phi}:H^2(\mathbb{U})\to H^2(\mathbb{U})$
are unitarily equivalent.
 \end{prop}

A combination of this result with (\ref{NormEstimate02}) gives us the following corollary.
  \begin{cor}\label{normestHtwo} If $T_{\psi,\phi}$ is a weighted composition operator on $H^2(\C_+)$, then
  $$
  \| T_{\psi,\phi}\|_{H^2(\C_+)} \le \|\Psi\|_\infty\sqrt{\frac {1+|\Phi(0)|}{1-|\Phi(0)|}}.
  $$
  \end{cor}

 One of the easiest compactness criteria for weighted composition operators is calculating their Hilbert--Schmidt norm, $\| \cdot \| _{\rm HS}$ and proving that it is finite. The formula used for that norm (see \cite{Valentin}), is:
  \[
  \| T_{\psi,\phi}\|_{\rm HS}=\sqrt{\frac{1}{4\pi }\int _{-\infty }^\infty \frac{|\psi (it)|^2}{{\textsf{ Re } }\phi (it)}\, dt}.
  \]
This principle will be used in the sequel.

 Also, we record for later use the fact that the unitary operator inducing the unitary equivalence in Proposition \ref{propValentin} is:  for $g\in H^2(\mathbb{U})$,
  \begin{equation}\label{Unitary}
Vg(w):=\frac{1}{1+w}g\left( \frac{w-1}{w+1}\right)\qquad  w\in \mathbb{C} _+.
  \end{equation}

 \subsection{The split}
 Our transfer operator is the sum of weighted composition operators
 $$
\L_p f(z) = p\sum_{n=1}^\infty  T_{\psi_n,\phi_n} f(z)+(1-p)\sum_{n=1}^\infty  T_{\widetilde\psi_n,\widetilde\phi_n} f(z),
 $$
 where for $n\ge 1$,
 $$
 \psi_n(z)=\widetilde \psi_n(z)=\frac {1}{(z+n)^2}, \text{ and } \phi_n(z)=\frac {1}{z+n},
 \quad\ \widetilde{\phi}_n= \frac {z+n-1}{z+n}.
 $$
  Similar to the previous section, we split $\L_p$ into sum of two operators:
 $$
 \B_p= (1-p) T_{\widetilde\psi_1,\widetilde\phi_1} f(z),\quad \A_p f(z)=\L_p f(z)-\B_pf(z).
 $$
 Moreover, let us also consider the split of $\A_p$
 $$
 \A^G_p f(z)= p\L_G f(z) = p\sum_{n=1}^\infty  T_{\psi_n,\phi_n} f(z),\quad
 \A^R_p f(z)= (1-p)\sum_{n=2}^\infty  T_{\widetilde\psi_n,\widetilde\phi_n} f(z).
 $$
 Here $\L_G$ is the operator as in \eqref{q:gaussrenyipf}.

  \begin{theorem}\label{compact} For any $p\in (0,1)$ the following three properties hold.
  \begin{enumerate}
  \item $\A_p^G=p\L_G:H^2(\C_+)\to H^2(\C_+)$ is compact.
  \item $\A_p^R:H^2(\C_+)\to H^2(\C_+)$ is nuclear.
  \item $J_p=(1-\B_p)^{-1}=\sum_{n=0}^\infty \B_p^n:H^2(\C_+)\to H^2(\C_+)$ is bounded.
  \end{enumerate}
  Therefore, $\widehat  \L_p = \A_p J_p=(\A_p^G+\A_p^R)J_p:H^2(\C_+)\to H^2(\C_+)$ is a compact operator.

 \end{theorem}

 \begin{cor} The unique invariant probability density $h_p$ on $[0,1]$ of the transfer operator $\L_p$ can be extended to an analytic function in $H^2(\C_+)$.
 \end{cor}

 \begin{proof}[Proof of Theorem \ref{compact}] We start with analysing the operator $\A_p^G$.
 We apply Corollary \ref{normestHtwo} to estimate the norms $\| T_{\psi_n,\phi_n}\|$ for  $n\ge 1$.

 For every $n\ge 1$, let $\Phi_n(z)=\gamma^{-1}\circ\phi_n\circ\gamma(z)$. Thus
 $$
 \Phi_n(0) = \frac {\phi_n(1)-1}{\phi_n(1)+1}=
 \frac
 {\frac {1}{n+1}-1}
 {\frac {1}{n+1}+1}=-\frac {n}{n+2},\text{ and } \sqrt{\frac {1+|\Phi_n(0)|}{1-|\Phi_n(0)|}}=\sqrt{n+1}.
 $$
 Moreover, if $w=\gamma(z)$, then $z=\gamma^{-1}(w)$,  and hence
 \[ \begin{split}
 \Psi_n(z)=\, & \psi_n (w) \frac{1-\frac {\phi_n(w)-1}{\phi_n(w)+1}}{1-\frac{w-1}{w+1}}= \psi_n (w) \frac {w+1}{\phi_n(w)+1}\\
 =\, & \frac {1}{(w+n)^2}\frac {(w+1)(w+n)}{w+n+1}=\frac {w+1}{(w+n)(w+n+1)}.
\end{split}\]
 Thus
 $$
 \sup_{|z|=1}|\Psi_n(z)|=\sup_{w=it, t\in\R} \left| \frac {w+1}{(w+n)(w+n+1)}\right|
 =\sup_{t\in\R} \sqrt{\frac{ (t^2+1)}{(t^2+n^2)(t^2+(n+1)^2)}} 
 =\frac 1{n(n+1)},$$
 and hence,
 $$\| T_{\psi_n,\phi_n}\| \le  \sup_{|z|=1}|\Psi_n(z)|\sqrt{\frac {1+|\Phi_n(0)|}{1-|\Phi_n(0)|}}= \frac 1{n\sqrt{n+1}}.
 $$
 Therefore $\A_p^Gf(z)=p\sum_{n=1}^\infty T_{\psi_n,\phi_n} f(z)$ is a sum of compact operators
 with
 $$
 \sum_{n=1}^\infty \| T_{\psi_n,\phi_n}\|<\infty,
 $$
 and hence $\A_p^G$ is compact. The reason why the operators $T_{\psi_n,\phi_n}$ are compact is that they are actually Hilbert--Schmidt. Indeed, one has that
 \[
 \| T_{\psi _n,\phi _n}\|_{\rm HS}=\sqrt{\frac{1}{4\pi }\int _{-\infty }^\infty \frac{|\psi _n(it)|^2}{{\textsf{ Re } }\phi _n(it)}\, dt}=\frac{1}{2n}<\infty \quad n=1,2,\dots
 \]
To check the above equality, note that
$$\aligned
\| T_{\psi _n,\phi _n}\|_{\rm HS}&=\sqrt{\frac{1}{4\pi }\int _{-\infty }^\infty \frac{|\psi _n(it)|^2}{{\textsf{ Re } }\phi _n(it)}\, dt}=
\sqrt{\frac{1}{4\pi }\int _{-\infty }^\infty \frac{\frac{1}{|n+it|^4}}{\textsf{ Re }\frac{1}{n+it}}\, dt}\\
&=\sqrt{\frac{1}{4\pi n}\int _{-\infty }^\infty \frac{dt}{n^2+t^2}}
=
\sqrt{\frac{1}{4\pi n^2}\arctan u|_{-\infty }^\infty} =\frac{1}{2n}.
\endaligned
$$

 We now turn to the proof of the second statement. We are going to show that  for any
 $n\ge 1$,
\begin{equation}\label{normestimtwo}
 \left\|T_{\widetilde \psi_n,\widetilde \phi_n}\right\|_1=\left\|T_{\widetilde \Psi_n,\widetilde \Phi_n}\right\|_1\le  \frac 2{n^{3/2}},
\end{equation}
 and hence,
 $$
 \| \A_p^R\|_1 \le (1-p)\sum_{n=2}^\infty \|T_{\widetilde \psi_n,\widetilde \phi_n}\|_1<\infty,
 $$
 and thus $\A_p^R$ is nuclear.

 Fix $n\ge 1$, then
 $$
 \widetilde \phi_n=1-\frac {1}{z+n}=\frac {z+n-1}{z+n}, \text{ and }
 \widetilde \Phi_n(z)=\gamma^{-1}\circ\widetilde \phi_n\circ\gamma(z)=-\frac {1-z}{(2n+1)-(2n-3)z}.
 $$
 Note that $\widetilde\Phi_n$ is a M\"obius transformation leaving the real line invariant. Since the real line is perpendicular to the unit circle, the circle will be transformed by $\widetilde\Phi_n$ into a circle perpendicular to the real line.
 Note also that $\widetilde\Phi_n(1)=0$ and $\widetilde\Phi_n(-1)=-1/(2n-1)$, hence $0$ and $-1/(2n-1)$ are antipodal points on that circle (i.e. the endpoints of a diameter). This makes it geometrically evident that
 \begin{equation}\label{NormInfinity}
\| \widetilde\Phi_n\| _\infty =\frac{1}{2n-1}\text{ and }
\frac{1}{1-\| \widetilde\Phi_n\| _\infty }=\frac{2n-1}{2n-2}.
\end{equation}
Next, we turn to the computation of $\| \widetilde{\Psi}_n\|_{H^2(\mathbb U)}$, where
$$\widetilde\Psi_n(z)=\frac{2(z-1)}{[(n-1)z-(n+1)][nz-(n+2)]}.
$$
By (\ref{Unitary}), one has that $\| \widetilde{\Psi}_n\|_{H^2(\mathbb U)}=\| V\widetilde{\Psi}_n\|_{H^2(\mathbb{C_+})}$. By a straightforward computation, one gets that
\[
V\widetilde{\Psi}_n(w)=-\frac{1}{(n+w)(n+1+w)},\qquad w\in \mathbb{C_+},
\]
and thus
\[
\aligned
\| V\widetilde{\Psi}_n\|_{H^2(\mathbb{C_+})}&=\sqrt{\frac{1}{\pi }\int_{-\infty }^\infty \frac{dt}{|n+it|^2|(n+1)+it|^2}}\\
&=\sqrt{\frac{1}{\pi (2n+1)}\int_{-\infty }^\infty \left( \frac{1}{n^2+t^2}-\frac{1}{(n+1)^2+t^2}\right)\, dt}\\
&=
\sqrt{\frac{1}{\pi (2n+1)}\left( \frac{1}{n}\arctan u\bigl|_{-\infty}^\infty- \frac{1}{n+1}\arctan v\bigl|_{-\infty}^\infty \right)}\\
&=
\left(\frac{ 1}{2n^3+3n^2+n}\right)^{\frac{1}{2}}\le \frac{1}{n^{3/2}},
\endaligned
\]
and hence (\ref{normestimtwo}) follows.

Let us now turn to the last statement. The operator
$$
\B f(z)=\frac 1{(1+z)^2} f\left( \frac {z}{1+z}\right)
$$
acts on $H^2(\C_+)$. We need to show that  the operator $
J_p=(1-(1-p)\B)^{-1} =\sum_{n=0}^\infty (1-p)^n\B^n$ is bounded.

Naturally,
$$
\left\| \sum_{n=0}^\infty (1-p)^n \B^n \right\|\le \sum_{n=0}^\infty (1-p)^n\left\|  \B^n \right\|.
$$
The second power satisfies
$$
\B^2f(z)=\frac {1}{(1+z)^2} \frac {1}{(1+\frac z{1+z})^2} f\left(\frac{ \frac {z}{1+z}}{
1+ \frac {z}{1+z}}\right)=\frac {1}{(1+2z)^2}f\left( \frac {z}{1+2z}\right)$$
and more generally, for any $n\ge 1$,
$$
\B^nf(z)= \frac {1}{(1+nz)^2}f\left( \frac {z}{1+nz}\right) =T_{\psi_n,\phi_n}f(z),
\quad
\phi_n(z) =\frac {z}{1+nz}, \quad \psi_n(z)=\frac {1}{(1+nz)^2}.
$$
By Proposition \ref{propValentin}, the corresponding unitarily equivalent weighted composition operator on $H^2(\mathbb U)$  is $T_{\Psi_n, \Phi_n}$  given by
$$
\Phi_n(z)=- \frac {n(z+1)}{(n-2)z+n+2},\quad \Psi_n(z)=-\frac {2(z-1)}{((n-2)z+n+2)((n-1)z+n+1)}.
$$
Hence $\| \Psi_n(z) \|_\infty=|\Psi_n(-1)|=|-1|=1$ and
$$
\| \B^n \|=\|T_{ \psi_n,\phi_n}\|=\|T_{\Psi_n,\Phi_n}\| \le \|\Psi_n\|_\infty\cdot \sqrt{
\frac {1+|\Phi_n(0)|}{1-|\Phi_n(0)|}}
\le 1\cdot \sqrt{ \frac {1+\frac n{n+2}}{1-\frac n{n+2}}}=\sqrt{ n+1},
$$
which allows one to conclude that
\[
\sum_{n=0}^\infty (1-p)^n\left\|  \B^n \right\|\le \sum_{n=0}^\infty (1-p)^n \sqrt{n+1}<\infty. \qedhere
\]
 \end{proof}

 \subsection{Real variable approach}
 The original approach of \cite{Bab, MR1} to the study of $\L_G$ on $H^2(\C_+)$ is based
 on the so-called real-variable theory of Hardy spaces. The Payley-Wiener theorem states
 that $f\in H^2(\C_+)$ if and only if there exists $g\in L^2(\R_+)$ such that
\begin{equation}\label{transform}
 f(z) =\int_{0}^{\infty} e^{-zt} g(t) dt,\quad \textsf{Re} z>0,
 \text{ and }\| f\|_{H^2(\C_+)}=  { {2\pi}} \| g\|_{L^2(\R_+)}.
\end{equation}
 Thus properties of operators acting on $H^2(\C_+)$ can equivalently be studied by considering the corresponding operators acting on $L^2(\R_+)$.
 A simple computation shows that if $f\in H^2(\C_+)$, satisfying \eqref{transform}, then for all $n\ge 1$,
\[ \begin{split}
 \frac{1}{(z+n)^2}f\left(\frac 1{z+n}\right)=\, &  \frac {1}{(z+n)^2}\int_0^\infty e^{-\frac t{z+n}} g(t)dt  =\int_0^\infty\sum_{k=0}^\infty \frac {(-t)^k}{k!}\frac 1{(z+n)^{k+2}} g(t) dt\\
 =\, & \int_0^\infty\sum_{k=0}^\infty \frac {(-t)^k}{k!} \left[\int_0^\infty e^{-zs}\frac {s^{k+1}e^{-ns}}{(k+1)!} ds\right]g(t) dt\\
 =\, & \int_{0}^\infty e^{-zs} e^{-ns} s  \left[\int_{0}^\infty \sum_{k=0}^\infty \frac{ (-1)^kt^ks^k}{k!(k+1)!}g(t)dt\right] ds\\
=\, & \int_{0}^\infty e^{-zs} e^{-ns} s \left[\int_{0}^\infty \frac {J_1(2\sqrt{st})}{\sqrt{st}} g(t)dt\right] ds,
 \end{split}\]
 and similarly,
\[ \begin{split}
 \frac{1}{(z+n)^2}f\left(1-\frac 1{z+n}\right)=\, & \frac {1}{(z+n)^2}\int_0^\infty e^{-t+\frac t{z+n}} g(t)dt =\int_0^\infty e^{-t}\sum_{k=0}^\infty \frac {t^k}{k!}\frac 1{(z+n)^{k+2}} g(t) dt\\
=\, & \int_{0}^\infty e^{-zs} e^{-(n-1)s} s \left[\int_{0}^\infty \frac {I_1(2\sqrt{st})}{\sqrt{st}} e^{-(s+t)}g(t)dt\right] ds,
 \end{split}\]
 where $J_1$ and $I_1$ are the Bessel and the modified Bessel functions
 of the first kind, respectively.
Therefore,
$$
\A_p^G f(z) =p\int_0^\infty e^{-zs} \frac s{e^s-1}\mathcal K_J g(s) ds,
\text{ and }
\A_p^R f(z)=(1-p)\int_0^\infty e^{-zs} \frac s{e^s-1}\mathcal K_I g(s) ds,
$$
 where
 $$
 \mathcal K_J g(s) =\int_{0}^\infty \frac {J_1(2\sqrt{st})}{\sqrt{st}} g(t)dt,
 \quad
  \mathcal K_I g(s) =\int_{0}^\infty \frac {I_1(2\sqrt{st})}{\sqrt{st}}e^{-(s+t)} g(t)dt
 $$
 are well-defined bounded linear operators on $L^2(\R_+)$.
 Moreover, since the kernels of these operators are symmetric,
 $ \mathcal K_J $ and $ \mathcal K_I$ are self-adjoint, and hence have real spectrum. Note however, that since
 $$
 \int_{0}^\infty \int_{0}^\infty\left| \frac {J_1(2\sqrt{st})}{\sqrt{st}}\right|^2 ds dt=+\infty,
 $$
 $\mathcal K_J$ is not Hilbert-Schmidt.

We note that  since $e^{-
\frac s2} \le \frac {s}{e^s-1}$ for all $s\ge 0$, for any $f\in H^2(\C_+)$, both $A_p^G f(z)$, $A_p^R f(z)$ are elements of the Hardy space $ H^2(\C_{-1/2})$,
where $\C_{-1/2}=\{z:\ \textsf{Re }z>-\frac 12\}$. In other words,
the operators $A_p^G f(z)$, $A_p^R f(z)$ do improve the `smoothness'.
However, singularities associated with $\B_p$ (c.f., expression for $\B_p^n$), result in that the eigenfunction is only in $H^2(\C_+)$,
and not in $H^2(\C_{-1/2})$, which was the case the transfer operator of the Gauss map.

The real-variable approach also gives a possibility of getting sharper
estimates of the norms of operators.
In the previous section we showed that $\A_p^R$ is nuclear ($\| \A_p^R\|_1<\infty$). Via the integral representation, we can derive sharper bounds: given the representation of the Bessel function $I_1$
as the series
$$\frac {I_1(2\sqrt{st})}{\sqrt{st}} e^{-(s+t)} = \sum_{n=0}^\infty\frac {t^n s^n }{n!(n+1)!} e^{-(t+s)},
$$
one has the following nuclear decomposition on $L^2(\R_+)$:
$$\aligned
\frac s{e^s-1} \mathcal K_I g(s) &=\sum_{n=0}^\infty \frac {s }{e^s-1} \frac {s^n e^{-s}}{(n+1)!}
\left\langle \frac {t^ne^{-t}}{n!}, g(t)\right\rangle_{L^2(\R_+)}=:\sum_{n=0}^\infty \xi_n(s) \left\langle \eta_n, g\right\rangle.\endaligned
$$
Similar to \cite{MR1}, one can easily show that
$$
\| \xi_n\|^2_{L^2(\R_+)} \le C \left(\frac 23\right)^{2n},\text{ and } \| \eta_n\|^2_{L^2(\R_+)} \le \frac {1}{\sqrt{n+1}},
$$
thus $\| \xi_n\|_{L^2(\R_+)}\cdot \| \eta_n\|_{L^2(\R_+)}  \le C\theta^n$ for
some $C>0$ and $\theta\in (0,1)$, and hence the operator
$$\mathcal  M_I g(s) = \frac {s}{e^s-1} \mathcal K_I g(s)
$$
is nuclear of order $0$.

 \section{Conclusions}
 The Gauss and Gauss-type continued fraction maps are classical examples in ergodic theory, and often serve as a play ground for the development and testing of various techniques.
 In the present paper we have applied
 such techniques to the study the transfer operator of the Gauss-R\'enyi random continued fractions
 map. We have shown that these techniques allow to conclude that this random map
 has an invariant density with roughly the same analytic properties, however,
 there is a clear 'loss' in smoothness due to the presence of an  indifferent fixed point.

\end{document}